\documentclass[12pt]{article}
\usepackage{mathrsfs}
\usepackage{amsmath,amssymb}

\usepackage{amsfonts}
\usepackage{latexsym}
\usepackage{amsthm}
\usepackage{pictex}

\usepackage[T2A]{fontenc}     
\usepackage[cp1251]{inputenc} 
\usepackage{graphicx}

\newcommand{\er}[1]{{\rm(\ref{#1})}}
\def\lb{\label}

\newtheorem{theorem}{Theorem}[section]
\newtheorem{definition}{Definition}
\newtheorem{lemma}{Lemma}[section]

\begin{document}

\def\a{\alpha} \def\cA{{\cal A}} \def\bA{{\bf A}}  \def\mA{{\mathscr A}}
\def\b{\beta}  \def\cB{{\cal B}} \def\bB{{\bf B}}  \def\mB{{\mathscr B}}
\def\g{\gamma} \def\cC{{\cal C}} \def\bC{{\bf C}}  \def\mC{{\mathscr C}}
\def\G{\Gamma} \def\cD{{\cal D}} \def\bD{{\bf D}}  \def\mD{{\mathscr D}}
\def\d{\delta} \def\cE{{\cal E}} \def\bE{{\bf E}}  \def\mE{{\mathscr E}}
\def\D{\Delta} \def\cF{{\cal F}} \def\bF{{\bf F}}  \def\mF{{\mathscr F}}
\def\c{\chi}   \def\cG{{\cal G}} \def\bG{{\bf G}}  \def\mG{{\mathscr G}}
\def\z{\zeta}  \def\cH{{\cal H}} \def\bH{{\bf H}}  \def\mH{{\mathscr H}}
\def\e{\eta}   \def\cI{{\cal I}} \def\bI{{\bf I}}  \def\mI{{\mathscr I}}
\def\p{\psi}   \def\cJ{{\cal J}} \def\bJ{{\bf J}}  \def\mJ{{\mathscr J}}
\def\vT{\Theta}\def\cK{{\cal K}} \def\bK{{\bf K}}  \def\mK{{\mathscr K}}
\def\k{\kappa} \def\cL{{\cal L}} \def\bL{{\bf L}}  \def\mL{{\mathscr L}}
\def\l{\lambda}\def\cM{{\cal M}} \def\bM{{\bf M}}  \def\mM{{\mathscr M}}
\def\L{\Lambda}\def\cN{{\cal N}} \def\bN{{\bf N}}  \def\mN{{\mathscr N}}
\def\m{\mu}    \def\cO{{\cal O}} \def\bO{{\bf O}}  \def\mO{{\mathscr O}}
\def\n{\nu}    \def\cP{{\cal P}} \def\bP{{\bf P}}  \def\mP{{\mathscr P}}
\def\r{\rho}   \def\cQ{{\cal Q}} \def\bQ{{\bf Q}}  \def\mQ{{\mathscr Q}}
\def\s{\sigma} \def\cR{{\cal R}} \def\bR{{\bf R}}  \def\mR{{\mathscr R}}
\def\S{\Sigma} \def\cS{{\cal S}} \def\bS{{\bf S}}  \def\mS{{\mathscr S}}
\def\t{\tau}   \def\cT{{\cal T}} \def\bT{{\bf T}}  \def\mT{{\mathscr T}}
\def\f{\phi}   \def\cU{{\cal U}} \def\bU{{\bf U}}  \def\mU{{\mathscr U}}
\def\F{\Phi}   \def\cV{{\cal V}} \def\bV{{\bf V}}  \def\mV{{\mathscr V}}
\def\P{\Psi}   \def\cW{{\cal W}} \def\bW{{\bf W}}  \def\mW{{\mathscr W}}
\def\o{\omega} \def\cX{{\cal X}} \def\bX{{\bf X}}  \def\mX{{\mathscr X}}
\def\x{\hat{\xi}}    \def\cY{{\cal Y}} \def\bY{{\bf Y}}  \def\mY{{\mathscr Y}}
\def\X{\hat{\Xi}}    \def\cZ{{\cal Z}} \def\bZ{{\bf Z}}  \def\mZ{{\mathscr Z}}
\def\O{\Omega}
\def\ve{\varepsilon}
\def\vt{\vartheta}
\def\vp{\varphi}
\def\vk{\varkappa}

\def\ba{\breve{a}}
\def\bb{\breve{b}}


\def\mmw{{\mathrm w}}
\def\mms{{\mathrm s}}
\def\mmg{{\mathrm g}}
\def\mmA{{\mathrm A}}
\def\mmJ{{\mathrm J}}


\newcommand{\gA}{\mathfrak{A}}
\newcommand{\gB}{\mathfrak{B}}
\newcommand{\gC}{\mathfrak{C}}
\newcommand{\gD}{\mathfrak{D}}
\newcommand{\gE}{\mathfrak{E}}
\newcommand{\gF}{\mathfrak{F}}
\newcommand{\gG}{\mathfrak{G}}
\newcommand{\gH}{\mathfrak{H}}
\newcommand{\gI}{\mathfrak{I}}
\newcommand{\gJ}{\mathfrak{J}}
\newcommand{\gK}{\mathfrak{K}}
\newcommand{\gL}{\mathfrak{L}}
\newcommand{\gM}{\mathfrak{M}}
\newcommand{\gN}{\mathfrak{N}}
\newcommand{\gO}{\mathfrak{O}}
\newcommand{\gP}{\mathfrak{P}}
\newcommand{\gR}{\mathfrak{R}}
\newcommand{\gS}{\mathfrak{S}}
\newcommand{\gT}{\mathfrak{T}}
\newcommand{\gU}{\mathfrak{U}}
\newcommand{\gV}{\mathfrak{V}}
\newcommand{\gW}{\mathfrak{W}}
\newcommand{\gX}{\mathfrak{X}}
\newcommand{\gY}{\mathfrak{Y}}
\newcommand{\gZ}{\mathfrak{Z}}

\newcommand{\gr}{\mathfrak{r}}
\newcommand{\gs}{\mathfrak{s}}


\def\mA{{\mathscr A}}
\def\mB{{\mathscr B}}
\def\mC{{\mathscr C}}
\def\mD{{\mathscr D}}
\def\mE{{\mathscr E}}
\def\mF{{\mathscr F}}
\def\mG{{\mathscr G}}
\def\mH{{\mathscr H}}
\def\mI{{\mathscr I}}
\def\mJ{{\mathscr J}}
\def\mK{{\mathscr K}}
\def\mL{{\mathscr L}}
\def\mM{{\mathscr M}}
\def\mN{{\mathscr N}}
\def\mO{{\mathscr O}}
\def\mP{{\mathscr P}}
\def\mQ{{\mathscr Q}}
\def\mR{{\mathscr R}}
\def\mS{{\mathscr S}}
\def\mT{{\mathscr T}}
\def\mU{{\mathscr U}}
\def\mV{{\mathscr V}}
\def\mW{{\mathscr W}}
\def\mX{{\mathscr X}}
\def\mY{{\mathscr Y}}
\def\mZ{{\mathscr Z}}

\def\Z{{\Bbb Z}}
\def\R{{\Bbb R}}
\def\C{{\Bbb C}}
\def\T{{\Bbb T}}
\def\N{{\Bbb N}}
\def\S{{\Bbb S}}
\def\H{{\Bbb H}}
\def\J{{\Bbb J}}
\def\dD{{\Bbb D}}
\def\W{{\Bbb W}}

\def\qqq{\qquad}
\def\qq{\quad}
\newcommand{\ma}{\begin{pmatrix}}
\newcommand{\am}{\end{pmatrix}}
\newcommand{\ca}{\begin{cases}}
\newcommand{\ac}{\end{cases}}
\let\ge\geqslant
\let\le\leqslant
\let\geq\geqslant
\let\leq\leqslant
\def\ma{\left(\begin{array}{cc}}
\def\am{\end{array}\right)}
\def\iint{\int\!\!\!\int}
\def\lt{\biggl}
\def\rt{\biggr}
\let\geq\geqslant
\let\leq\leqslant
\def\[{\begin{equation}}
\def\]{\end{equation}}
\def\wh{\widehat}
\def\wt{\widetilde}
\def\pa{\partial}
\def\sm{\setminus}
\def\es{\emptyset}
\def\no{\noindent}
\def\ol{\overline}
\def\iy{\infty}
\def\ev{\equiv}
\def\/{\over}
\def\ts{\times}
\def\os{\oplus}
\def\ss{\subset}
\def\h{\hat}
\def\Re{\mathop{\rm Re}\nolimits}
\def\Im{\mathop{\rm Im}\nolimits}
\def\supp{\mathop{\rm supp}\nolimits}
\def\sign{\mathop{\rm sign}\nolimits}
\def\Ran{\mathop{\rm Ran}\nolimits}
\def\Ker{\mathop{\rm Ker}\nolimits}
\def\Tr{\mathop{\rm Tr}\nolimits}
\def\const{\mathop{\rm const}\nolimits}
\def\dist{\mathop{\rm dist}\nolimits}
\def\diag{\mathop{\rm diag}\nolimits}
\def\Wr{\mathop{\rm Wr}\nolimits}
\def\BBox{\hspace{1mm}\vrule height6pt width5.5pt depth0pt \hspace{6pt}}

\def\Diag{\mathop{\rm Diag}\nolimits}


\def\ve{\varepsilon}   \def\vt{\vartheta}    \def\vp{\varphi}    \def\vk{\varkappa}

\def\Z{{\mathbb Z}}    \def\R{{\mathbb R}}   \def\C{{\mathbb C}}    \def\K{{\mathbb K}}
\def\T{{\mathbb T}}    \def\N{{\mathbb N}}   \def\dD{{\mathbb D}}    \def\S{{\mathbb S}}


\def\Twelve{
\font\Tenmsa=msam10 scaled 1200 \font\Sevenmsa=msam7 scaled 1200
\font\Fivemsa=msam5 scaled 1200 \textfont\msbfam=\Tenmsb
\scriptfont\msbfam=\Sevenmsb \scriptscriptfont\msbfam=\Fivemsb

\font\Teneufm=eufm10 scaled 1200 \font\Seveneufm=eufm7 scaled 1200
\font\Fiveeufm=eufm5 scaled 1200
\textfont\eufmfam=\Teneufm \scriptfont\eufmfam=\Seveneufm
\scriptscriptfont\eufmfam=\Fiveeufm}

\def\Ten{
\textfont\msafam=\tenmsa \scriptfont\msafam=\sevenmsa
\scriptscriptfont\msafam=\fivemsa

\textfont\msbfam=\tenmsb \scriptfont\msbfam=\sevenmsb
\scriptscriptfont\msbfam=\fivemsb

\textfont\eufmfam=\teneufm \scriptfont\eufmfam=\seveneufm
\scriptscriptfont\eufmfam=\fiveeufm}

\title {Periodic Jacobi operator with  finitely supported perturbations: the inverse resonance  problem.}

\author{ Alexei Iantchenko
\begin{footnote}
{Malm{\"o} H{\"o}gskola, email: ai@mah.se }
\end{footnote}
\and
Evgeny Korotyaev
\begin{footnote}
{Saint-Petersburg University, e-mail: korotyaev@gmail.com}
\end{footnote}
}

\maketitle

\begin{abstract}
 \no We consider a periodic Jacobi operator $H$ with finitely supported perturbations  on $\Z.$
  We solve the inverse resonance problem: we prove that the mapping from finitely supported
 perturbations to the scattering data, the inverse of the transmission coefficient and the  Jost function on the right half-axis,  is one-to-one and onto. We consider the problem of reconstruction of the scattering data from all  eigenvalues, resonances and the set of zeros of $R_-(\l)+1,$   where $R_-$ is the reflection coefficient.
\end{abstract}


\section{Introduction.}
\setcounter{equation}{0}

We consider a  Jacobi operator $H=H^0+V$ on the lattice
$\Z=\{\ldots,-2,-1,0,1,2,\ldots\}$. Here the unperturbed operator $H^0$ is a periodic
Jacobi operator given by
\[
\lb{J0}
 (H^0y)_n={a}_{n-1}^0y_{n-1}+{a}_{n}^0y_{n+1}+{b}_n^0y_n,
 \]
where  $y=(y_n)_1^\iy\in \ell^2=\ell^2(\Z)$ and   the $q-$periodic
coefficients $a_n^0, b_n^0\in \R$ satisfy
  \[
 \lb{1e}
{a}^0_n={a}^0_{n+q}>0,\qq b_n^0={b}^0_{n+q},\qq n\in
\Z,\qqq \prod_{j=1}^qa_j^0=1,\qq q\ge 2.
\]
We fix a positive integer $p.$ The perturbation operator $V$ is the finitely supported Jacobi
operator given by
\[
\lb{V} (Vy)_n=\ca {u}_{n-1}y_{n-1}+u_{n}y_{n+1}+v_ny_n, & {\rm
if}\qq  1\le n\le p,\\
u_py_p,& {\rm if} \qq n=p+1,\\
u_0y_1+v_0y_0,& {\rm if} \qq n=0,\\
0, \qqq & {\rm if} \qq n\leq -1\,\,\mbox{or}\,\,n\ge p+2,\qq p\ge 1.\ac
\]

 We parameterize $V$ by the vector $(u,v)\in \R^{2p}$ and let $(u,v)$
belong to the class $\cV_\n$ given by
\begin{align}
\label{class}
 \cV_\nu=&\lt\{ (u,v)\in \R^{2p}:\,\,
a^0_n+u_n>0,\,\, n=0,...,p, \,\, u_p\neq0,\,\,v_0\neq 0\rt\}\,\, \mbox{if}\,\, \n=2p,
\\
\label{class2} \cV_\nu=&\lt\{ ({u} ,{v} )\in
\R^{2p}:\,\,{a}^0_n+{u} _n >0 ,\,\, n=0, ...,p,\,\,{u} _p= 0,\,\,v_0\neq 0,\,\, {v} _p\neq
 0,\,\,\rt\}\,\,\mbox{if}\,\, \n=2p-1.
\end{align}

We rewrite $H$ in the form
\[
\lb{pert}
 (H y)_n={a}_{n-1}y_{n-1}+{a}_{n}y_{n+1}+{b}_ny_n
 \]
with the coefficients $a_n, b_n$ given by
\[
\lb{ab}{a}_n=\ca {a}^0_n+{u}_n>0& {\rm if} \qq 0\le n\le p,\\
                           {a}^0_n & {\rm if} \qq n\le -1\,\,\mbox{or}\,\, n\ge p+1,
         \ac
              \qq {b}_n=\ca  b_n^0+{v}_n & {\rm if} \qq 0\le n\le p, \\
            b_n^0   & {\rm if} \qq n\le -1\,\,\mbox{or}\,\, n\ge p+1.  \ac
\]
 The  corresponding  Jacobi matrices have the forms
\[
\lb{Hm}
H^0=\left(\begin{array}{cccccccc}
  ...& ...  & ...  & ...     &...         & ...        &... \\
 ...&a_0^0  & b_1^0  & a_1^0     &0         & 0        &... \\
 ...&0  & a_1^0& b_2^0       &a_2^0& 0        &... \\
 ...&0  & 0  &a_2^0& b_3^0         & b_3^0     &... \\
 ...&0  & 0  & 0        &a_3^0      & b_4^0       &... \\
 ...&0  & 0  & 0        &0         &a_4^0&... \\
 ...&...  & ...& ...      &...       &...       &... \\
 \end{array}\right),\qqq
 H=\left(\begin{array}{cccccccc}
  ...& ...  & ...  & ...     &...         & ...        &... \\
 ...&a_0  & b_1  & a_1     &0         & 0        &... \\
 ...&0  & a_1& b_2       &a_2& 0        &... \\
 ...&0  & 0  &a_2& b_3         & b_3     &... \\
 ...&0  & 0  & 0        &a_3      & b_4       &... \\
 ...&0  & 0  & 0        &0         &a_4&... \\
 ...&...  & ...& ...      &...       &...       &... \\
 \end{array}\right).
\]

For $a_n=1,$ $b_n^0=0,$ $n\in \Z,$ the operator $H$ is the finite difference Schr{\"o}dinger operator with finitely supported potential.

 A lot of papers is devoted to the direct and inverse resonance problems for the Schr{\"o}dinger operator $-\frac{d^2}{dx^2}+q(x)$ on the line $\R$  with compactly supported perturbation (see \cite{S}, \cite{Fr}, \cite{Z}, \cite{K3}  and references given there).
Zworski [Z]
obtained the first results about the distribution of resonances for
the Schr\"odinger operator with compactly supported potentials on
the real line. One of the present authors obtained the uniqueness,
the recovery and the characterization of the $S$-matrix for the
Schr\"odinger operator with a compactly supported potential on the
real line \cite{K3}, see also
\cite{Z1}, \cite{BKW} concerning the uniqueness.

The problem of resonances for the Schr{\"o}dinger with periodic plus compactly supported potential
$-\frac{d^2}{dx^2}+p(x)+q(x)$  is much less studied: \cite{F1}, \cite{KM}, \cite{K1}.
The
following results were obtained in \cite{K1}: 1) the
distribution of resonances in the disk with large radius is
determined, 2) some inverse resonance problem, 3) the existence of a
logarithmic resonance-free region near the real axis.
The inverse resonance problem is not yet solved.

Finite-difference Schr{\"o}dinger and  Jacobi operators
express many similar features.
Spectral and scattering properties of infinite Jacobi matrices are much studied (see \cite{Mo}, \cite{DS1}, \cite{DS2} and references given there).
The inverse problem was solved for periodic Jacobi operators: \cite{P}.

 The inverse scattering problem  for asymptotically periodic
 coefficients was solved by  Khanmamedov: \cite{Kh1}
 (note that the russian versions were dated much earlier) and Egorova, Michor and Teschl \cite{EMT}
 (in the case of quasi-periodic background).

The resonance problems are less studied (see M.Marletta and R.Weikard \cite{MW}).
The inverse resonances problem was recently solved in the case of
constant background \cite{K2}.

In \cite{IK1} we consider the direct resonance problem in the case of periodic background. We describe the spectral and scattering properties of $H.$ Moreover,
in  the special case $ u_n\equiv 0$ we   obtain the asymptotics of the spectrum in the limit of small perturbations $ V.$ In Theorem \ref{ThIK1}  below we summarize some results obtained in \cite{IK1}.

In \cite{IK2} we consider the zigzag half-nanotubes
(tight-binding approximation) in a uniform magnetic field which is described by the magnetic Schr\"odinger operator with a periodic potential plus a finitely supported perturbation on the half-lattice.
We describe all eigenvalues and resonances of this operator, and their dependence on the magnetic field.

In \cite{IK3} we consider a periodic Jacobi operator  with finitely supported perturbations  on the half-lattice.
 We describe all eigenvalues and resonances,  and give
their properties. We solve the inverse resonance problem: we prove that the mapping from finitely supported
 perturbations to the Jost functions is one-to-one and onto, we show how the Jost functions can be reconstructed from all  eigenvalues, resonances and from the set of zeros of $S(\l)-1,$ where $S(\l)$ is the scattering matrix.

In the present paper we extend the methods from \cite{IK3} to the inverse resonance problem on the lattice $\Z.$ In one aspect the inverse scattering problem for the  perturbed operator $H$ on $\Z$
  is simpler then on the half-lattice:   even the unperturbed periodic Jacobi operator on the half-lattice has bound and antibound states. But technically the inverse problem on the lattice is more involved as we need to reconstruct two analytic functions on the two-sheeted Riemann surface: the numerator and denominator of the reflection coefficient $R_-,$ instead of one as in \cite{IK3}, the Jost function $f_0^+.$

Now we pass to the description of the spectral and stattering properties of $H,$ recalling some results from \cite{IK1}, and formulate our main results.

The
  spectrum of
$H^0$ on $\ell^2(\Z)$ is absolutely continuous and consists
of $q$ zones $\sigma_j$
separated by the gaps $\gamma_j$ given by
\begin{align}
 &\s_j=[\l_{j-1}^+,\l_j^-],\qq j=1,\ldots,q,\qq
\g_{j}=(\l^-_{j},\l^+_j),\qqq j=1,\ldots,q-1,\nonumber\\
&\l_0^+<\l_1^-\le \l_1^+<.. ...<\l^-_{q-1}\le \l^+_{q-1}
<\l^-_{q}. \lb{ends}
\end{align}

 We denote $\gamma_0=(-\infty, \l_0^+)$ and $\gamma_q=( \l_q^+,+\infty)$ the infinite gaps.

Let $\vp=(\vp_n(\l))_1^\iy$ and $ \vt=(\vt_n(\l))_1^\iy$ be
fundamental solutions for the equation
\[
\lb{1e2} a_{n-1}^0 y_{n-1}+{a}_{n}^0y_{n+1}+{b}_n^0y_n=\l y_n,\qqq
\l\in\C,
\]
satisfying the conditions $\vt_0=\vp_1=1$ and $\vt_1=\vp_0=0$. Here and
below $a_0^0=a_q^0$.  Introduce the  Lyapunov function  $\D$ by
\[
\lb{LF} \D={\vp_{q+1}+\vt_{q}\/2}.
\]
It is known that $\D(\l)$ is a polynomial of degree $q$ and
$\l_j^\pm, j=1,...,q,$ are the zeros of the polynomial $\D^2(\l)-1$
of degree $2q$. Note that $\D(\l_j^\pm)=(-1)^{q-j}$. In each ``gap''
$[\l_j^-, \l_j^+]$ there is one simple zero of polynomials $\vp_q,
\dot \D, \vt_{q+1}.$ Here and below $\dot f$ denotes the derivative
of $f=f(\l)$ with respect to $\l:$ $\dot f\equiv\partial_\l f\equiv
f'(\l).$

Let $\G$ denote the complex plane cut along the segments $\sigma_j$ (\ref{ends}):
 $\G=\C\setminus\s_{\rm ac}(H^0).$ Now we introduce the two-sheeted Riemann
surface $\L$  of $\sqrt{1-\D^2(\l)}$  by joining the upper and lower rims of
two copies of the cut plane $\G$ in
the usual (crosswise) way.  We identify the first (physical) sheet $\L_1$ with $\G$
and the second sheet we denote by $\L_2$.

 Let $\,\,\wt{}\,\,$ denote  the natural projection from
 $\L$  into the complex plane:
 \begin{equation}
 \lb{projection}
\l\in\L,\qq\l\to \wt\l\in\C.
\end{equation}
By identification of $\G=\C\setminus\s_{\rm ac}(H^0)$ with $\L_1,$ map $\,\,\wt{}\,\,$ can be also considered to be projection from $\L$ into the physical sheet $\L_1.$

The $j-$th gap on the first physical sheet $\Lambda_1$ we will denote by $\g_j^+$ and the same gap but on
the second nonphysical sheet $\Lambda_2$  we will denote by $\g_j^-$ and let $\g^{\rm c}_j$  be the union of $\overline{\g^+_j}$ and $\overline{\g^-_j}$:
\[
\lb{union}
\g_j^{\rm c}=\overline{\g^+_j}\cup\overline{\g^-_j}.
\]
Define the function $\O(\l)= \sqrt{1-\D^2(\l)}, \l\in \L,$ by
\begin{equation}\lb{branch}
\O(\l)<0\qq\mbox{for}\qq\l\in (\l_{q-1}^+, \l_{q}^-)\ss \L_1.
\end{equation}
 Introduce the
Bloch functions $\psi_n^\pm$ and the Titchmarch-Weyl functions
$m_\pm$ on $\L$ by
\begin{align}
 &\psi_n^\pm (\l)=\vt_n(\l)+m_\pm(\l)\vp_n(\l),\lb{Jost}\\
&m_\pm(\l)= \frac{\f(\l)\pm i\O(\l)}{\vp_q},\,\qqq
\phi=\frac{\vp_{q+1}-\vt_q}{2},\,\,\l\in\L_1. \lb{TitchWeyl}
\end{align}
 The projection of all singularities of $m_\pm$ to the complex plane coincides with the set of zeros
$\{\mu_j\}_{j=1}^{q-1}$ of polynomial $\vp_q$. Recall that $\vt_n, \vp_n, \f$ are polynomials.
Recall that any polynomial  $P(\l)$  gives rise to a function $P(\l)=P(\wt\l)$  on
the  Riemann surface $\L$ of $\sqrt{1-\D^2(\l)}.$

The perturbation $V$ satisfying (\ref{V}) does not change the absolutely continuous spectrum:
\[
 \lb{acspectr}\s_{\rm ac}(H)=\s_{\rm ac}(H^0)=
 \bigcup_{n=1}^q[\l_{n-1}^+,\l_n^-].
\]
The spectrum of $H$ consists of an absolutely continuous part
$\s_{\rm ac}(H)=\s_{\rm ac}(H^0)$ plus a finite number of simple eigenvalues in
each non-empty gap $\g_n, n=0,...,q$.

Introduce the function $$\alpha(\l) = C\det ((H-\l)(H^0-\l)^{-1})=C\det
(I+(H-H^0)(H^0-\l)^{-1}),\qq C=\prod_{j=0}^p\frac{a^0_j}{ a_j},$$
 which is meromorphic on $\L$, see [F1].
Recall that $T=1/\alpha$ is the transmission coefficient in the
$S-$matrix for the pair $H,H^0$ (see Section \ref{s-pert}).   If $\alpha$ has some
poles, then they coincide with some $\l_k^\pm .$ It is well known
that if $\alpha(\l) =0$ for some zero $\l\in\L_1,$ then $\l$ is
an eigenvalue of $H$  and $\l\in\cup\gamma_k^+.$ Note that there are
no eigenvalues on the spectrum $\sigma_{\rm
ac}(H^0)\subset\L_1$ since $|\alpha (\l)|\geq 1$ on
$\sigma_{\rm ac} (H^0).$

We define the functions $A$, $J$  by
$$J(\l)=2\Omega(\l+i0)\Im\alpha(\l+i0),\qq A(\l)=\Re\alpha (\l+i0)
-1,\qq\mbox{for}\,\,\l\in\sigma\, (H^0)\subset \L_1.
$$
These functions were introduced for the  Schr{\"o}dinger operator on $\R$ with periodic plus compactly supported potentials  by the second author in \cite{K1}.
   We  show that $A, J$ are polynomials on $\C$ and they are real on the real
line. Instead of the function $\alpha$ we consider the modified
function $\hat{w} =2i\Omega\alpha$ on $\Lambda.$ We show that $\hat{w}$
satisfies
\begin{equation}\lb{1.3}\hat{w} =2i\Omega\alpha =2i\Omega (1+A)
-J\qq\mbox{on}\,\,\Lambda .
\end{equation}

Recall that $\Omega$
 is analytic on $\Lambda$ and
$\Omega = 0$ for some $\l\in\Lambda$ iff $\l =\l_k^-$ or $\l
=\l_k^+$ for some $k\geq
 0.$ Then the function $\hat{w}$ is analytic on $\Lambda$ and has branch points $\l_n^\pm$ if $\gamma_n\neq \emptyset.$  The zeros of $\hat{w}$  are the eigenvalues and
the resonances.  Define the set
$$\Lambda_0=\{\l\in\Lambda:\,\,\l
 =\l_k^+\in\L_1\,\,\mbox{and}\,\, \l
 =\l_k^+\in\L_2,\,\, \gamma_k=\emptyset\}\subset\Lambda.$$
  In fact with each $\gamma_k =\emptyset$ we associate two points
  $\l_k^+\in\L_1$ and $\l_k^+\in\L_2$ from the
set $\Lambda_0.$  If each gap of $H^0$ is not empty, then
$\Lambda_0=\emptyset.$

\begin{definition}\lb{states} Each zero of $\hat{w}$ on $\Lambda\setminus\Lambda_0$
is a state of $H.$\\
1) A state $\l\in\L_1$ is a bound state.\\
2) A state $\l\in\L_2$ is a resonance.\\
3) A state $\l=\l_k^\pm,$ $k=1,\ldots,q$ is a virtual state.\\
A resonance $\l\in \cup\gamma_k^-\subset\L_2$ is an anti-bound
state.
\end{definition}
It is known that the gaps $\gamma_k=\emptyset$ do not give
contribution to the states. Recall that $S-$matrix for $H,H^0$ is
meromorphic on $\Lambda,$ but it is analytic at the points from
$\Lambda_0$ (see [F1]). Roughly speaking there is no difference
between the points from $\Lambda_0$  and other points inside the
spectrum of $H^0.$


In accordance with the continuous case \cite{K1} we define the
important function
\[\lb{FF}
 \cF(\l)=\hat{w}(\l)\hat{w}^*(\l),\qqq \l\in \L_1,
\]
where  we put $f^*(\l):=\overline{f(\overline{\l})}.$
For the perturbation $V$ with $(u,v)\in\cV_\nu$ we define the
constants
\[
 \lb{c2_2p-1}
 c_3=c_1c_2,\qqq
 c_1={1\/\prod_{0}^p{a}_j},\qqq \qqq
 c_2=\ca c_1u_p(a^0_p+a_p) &\qq \mbox{if}\qqq \nu=2p,\\
c_1(a^0_p)^2{v}_p          & \qq\mbox{if}\qqq\nu=2p-1. \ac
\]
 The distribution of the states is summarized in
the following theorem (see \cite{IK1}).

\begin{theorem}\lb{ThIK1} Let  the  Jacobi operator $H=H^0+V$
satisfy \er{J0}--\er{V}. Suppose $(u,v)\in \cV_\n$, where $\n\in \{2p,2p-1\}$. Then $\hat{w}$ satisfies (\ref{1.3}) and the following facts hold true.\\
 i) The function $\cF(\l)=\hat{w}(\l)\hat{w}^*(\l),$ $\l\in \L_1,$ is a real polynomial. Each zero
of $\cF$ is the projection of a state of $H$ on the first sheet. There
are no other zeros. The multiplicity of a bound state and a resonance is the multiplicity of its projection as a zero of $\cF.$  All bound states are simple.
The virtual state at $\l_j^\pm,$ $\gamma_j\neq \emptyset,$ $j=1,\ldots, q-1,$ is a simple zero of $\cF.$ Moreover, $\cF$ satisfies
\[
\lb{tn} \cF(\l)=-\l^{\k}(c_3v_0+\cO(\l^{-1})),\qqq  \k=\n+2q-1,\qqq
\l\to \iy,
\]
here $\k$ is the total number of states (counted with
multiplicities).\\
ii) There exists an even number of states (counted with multiplicities) on each set $\gamma_j^{\rm c}\neq\emptyset,$ $j=1,\ldots,q-1,$ where $\gamma_j^{\rm c}$ is a union of the physical gap $\overline{\gamma^+_j}$ and non-physical gap $\overline{\gamma^-_j}$ (see (\ref{union})).\\
iii) Let $\l_1\in\gamma_j^+$ be a bound state for some $j=0,\ldots,q,$
i.e. $\hat{w}(\l_1)=0.$ Let $\l_2\in\gamma_j^-\subset\Lambda_2$ be the same number but on the second sheet $\Lambda_2.$ Then $\l_2\in\gamma_j^-$ is not an antibound state, i.e. $\hat{w}(\l_2)\neq 0.$
\end{theorem}

Let $f_n^+$ denote the Jost solution for the equation $Hy=\l y$ satisfying
$f_n^+=\psi_n^+$ for $n\geq p+1$ (see (\ref{peq1}), (\ref{defJostline})) and $f_0^+$ is called the Jost function.
We prove that the operator $H$ is uniquely determined by the pair $(\hat{w},f_0^+).$ In the following definition we describe the class of functions with characteristic properties of $(\hat{w},f_0^+).$

\begin{definition}\lb{class}
For $\nu\in\N,$ let $\gC_\nu$    denote  the class of  pairs of functions  $(\mmw, f)$  on $\Lambda:\\$
$\mmw$ is entire function  of the form
\begin{align*}
&   \mmw=2i\Omega(1+{\mathrm A})-{\mathrm J},\\
&  \mmw=\ca \frac{A_p}{c_1}\l^q\left(1+{\mathcal O}(\l^{-1})\right) \qqq &{\rm if}
\qqq   \l\in\L_1\\
 -\frac{v_0}{A_p} c_2\l^{\nu+q-1}\left(1+{\mathcal O}(\l^{-1})\right) \qqq &{\rm
if}\qqq \l\in\L_2 \ac \qqq{\rm as}\qq \l\to \iy,
\end{align*}
where   $\mmA,$
$\mmJ$ are real polynomials (with real coefficients) of the orders $\nu -1$ and $\nu +q-1$ respectively;\\
 $f$ is meromorphic  function of the form
\begin{align*}
& f=P_1+\frac{\phi}{\vp_q}P_2+ i\frac{\Omega(\l)}{\vp_q}P_2, \\
& f(\l)=\ca c_1A_p+{\cO}(\l^{-1}) \qqq &{\rm if}
\qqq   \l\in\L_1\\
 -\frac{c_2 }{A_p}\l^{\nu}+\cO(\l^{\nu-1}) \qqq &{\rm
if}\qqq \l\in\L_2 \ac \qqq{\rm as}\qq \l\to \iy,
\end{align*}
where
 $P_1$ and $P_2$ are real polynomials  of
the orders $\nu -2$ and $\nu -1$ respectively.

Here $A_p$ is given by  \[
\lb{defAp}
A_p=\prod_{n=0}^pa_n^0.
\] The real constants  $c_1,$ $c_2,$ $v_0$ satisfy:  $c_1>0,$ $c_2\ne 0,$ $v_0\ne 0.$

We denote
$\sigma_{\rm st}\subset\L$  the set of all zeros of $\mmw.$ Denote  $\sigma_{\rm bs}=\sigma_{\rm st}\cap\L_1,$ $\sigma_{\rm st}=\{\rho_k\}_{k=1}^N,$  and let $$\cF =\mmw \mmw^*=4(1-\Delta^2)(1+\mmA)^2+\mmJ^2,\qqq {\mathrm s}=\frac{2i\Omega}{a_0^0}f- \mmw.$$
Suppose that the following properties are satisfied:\\
\no 1) $\sigma_{\rm bs}\in \cup_0^q\g_j^+$\\
\no 2) The function $\cF$ has even
number of zeros on each interval $[\l_j^-,\l_j^+],$ $j=1,\ldots,q-1,$
and $\cF$ has only simple zeros
at $\l_j^\pm,$ $\gamma_j\neq\emptyset,$ and at $\rho_k\in\sigma_{\rm bs},$ $k=1,\ldots, N.$\\
\no 3) For $k=1,\ldots,N,$
\begin{equation}\lb{gamma_plus}
 \mmg_k=-\frac{D^-}{D^+}\frac{2i\Omega(\rho_k)}{\mms\mmw'(\rho_k)}>0.
 \end{equation} \\
\no  4)
If $\vp_q(\l_j^\pm)=0,$ for some  $j=0,\ldots, q,$ then $\mms(\l_j^\pm)=\mmw(\l_j^\pm).$

\end{definition}

 If $(u,v)\in \cV_\nu,$ then  Lemma \ref{l_wf_in_class} states that
  $(\hat{w},f_0^+)\in\gC_\nu$ with $f=f_0^+,$ $P_1=\vt_0^+,$  $P_2=\vp_0^+,$   $\mmw=\hat{w},$ $\mmA=A$ and $\mmJ=J.$ Moreover, $\mms=\hat{s}:=\vp_q s/a_0^0$ with $s$ defined in (\ref{def_of_sw}) and  $\mmg_k=\gamma_{+,k}$ is the norming constant defined in (\ref{gamma_plus}).

Now we construct the mapping $\gF: \cV_\nu\to \gC_\nu,$
$\nu\in\{2p-1, 2p\},$ by the rule:
\[
 ({u} ,{v} )\to (\hat{w},f_0^+),
 \]
  i.e.  to each
$({u} ,{v} )\in \cV_\nu$ we associate $(\hat{w},f_0^+)\in \gC_\nu$.

Our main  result is formulated in the following theorem.
\begin{theorem}
\lb{Th1}
Let $\k=\nu+2q-1,$ $\nu\in\{2p-1,2p\},$ and suppose $\k\geq 2q+1.$ Then the mapping $\gF:\,\,\cV_\nu\mapsto \gC_\nu$ given by $\gF(V)=(\hat{w},f_0^+)$ is one-to-one and onto. Moreover, the reconstruction algorithm is specified.
\end{theorem}

In Theorem \er{Th1} we solve the inverse problem for the mapping
$\gF.$ The solution is  divided into the following three parts.
\begin{enumerate}
\item Uniqueness. Do the  pair of functions  $(\hat{w},f_0^+)\in\gC_\nu$ determine uniquely
 $({u} ,{v} )\in\cV_\nu$?
\item Reconstruction. Give an algorithm for recovering $({u} ,{v} )$ from
$(\hat{w},f_0^+)\in \gC_\k$ only.
\item Characterization.  Give necessary and sufficient conditions for
$(\mmw,f)$ to be the Jost function and $\hat{w}$ for some perturbation $({u} ,{v} )\in\cV_\nu.$
    \end{enumerate}

Using the polynomial interpolation as in \cite{IK3}, Theorems 1.4 and 1.5, we get that  operator $H$ can be reconstructed from $H^0,$ $\s_{\rm st}(H),$ $Zeros\,(R_-+1)$ and two polynomials $A,$ $\vp_0^+$ which follows from the following theorem.
\begin{theorem}\lb{zerosR+1} Suppose that the functions $\hat{w}$  and $R_-+1$ have only simple zeros,  disjoint from the end-points $\l_k^\pm,$ $k=0,\ldots,q,$  and from the points $\mu_j,$ $j=1,\ldots,q-1,$ such that $\vp_q (\wt\mu_j)=0$ (the Dirichlet eigenvalues).  Then the pair of functions $(\hat{w},f_0^+)$ is uniquely determined by the bound states and resonances of $H,$ the set of zeros of the function $R_-+1,$    the polynomials  $A,$  $\vp_0^+,$ and the constants $c_3,$ $v_0$ (see (\ref{c2_2p-1})).
\end{theorem}
 Note that from the assumptions of Theorem \ref{zerosR+1} it follows that the set of zeros of $\hat{w}$ coincide with $\sigma_{\rm st}(H)= \sigma_{\rm bs}(H)\cup\sigma_{\rm r}(H)$ (no virtual states are present) and
   coincide with the zeros of $\alpha=1/T.$ Moreover, the zeros of $R_-+1$ coincide with the zeros of the Jost function $f_0^+$ (see the proof Theorem \ref{zerosR+1}).

{\bf Plan of the paper.}
 In Part \ref{s-prel} we recall the  construction of the quasi-momentum map and the associated Riemann surface (Section \ref{ss-cut}), and
 consider the scattering problem by finitely supported perturbations (Section \ref{s-scattering}).

In Part \ref{s-first}  we prove the Theorems \ref{Th1} and \ref{zerosR+1}.
In Section \ref{s-pert} we summarize the characteristic properties of the scattering data.
In Section \ref{s-invscat} we summarize  the inverse scattering results from  \cite{Kh1} and \cite{EMT} in the form suitable for us.
Finally in Sections \ref{s-Th1} and \ref{s-zerosR+1} we prove Theorems \ref{Th1}and \ref{zerosR+1}.

\section{Preliminaries.}\lb{s-prel}

\subsection{Quasi-momentum map and Riemann surface $\cZ.$}\lb{ss-cut}

 In this section we recall the construction of  the conformal mapping  of the Riemann surface onto
 the plan with ``radial slits'' $\cZ,$ given in \cite{IK1}.    Our definition corrects
 the similar construction  in \cite{BE} and \cite{EMT}, where there
 was a mistake.

We suppose that all  gaps are open:
$\lambda_j^-<\lambda_j^+,$ $j=1,\ldots,q-1$.

Introduce a domain $\C\sm\cup_0^{q}\ol\g_j$ and a quasi-momentum domain
$\K$ by
$$
\K=\{\vk\in \C:-\pi\le\Re \vk\le 0\}\sm\cup_1^{q-1}\ol\G_j,\ \ \G_j =\rt(
-{\pi j+ih_j\/q}, -{\pi j-ih_j\/q}\rt).
$$
Here $h_j\ge 0$ is defined by the equation $\cosh
h_j=(-1)^{j-q}\D(\alpha_j)$ and $\alpha_j$ is a zero of $\D'(\l)$ in
the ``gap'' $[\l_j^-,\l_j^+]$. For each periodic Jacobi operator
there exists a unique conformal mapping $\vk:\C\sm\cup_0^{q}\ol\g_j\to
\K$ such that the following identities and asymptotics hold true:
\[
\lb{ca}
  \cos q\vk(\l)=\D (\l),   \ \ \ \l\in \C\sm\cup_0^{q}\ol\g_j,
  \ \ \ \ \ {\rm and} \ \ \
  \ \ \ \vk(it)\to \pm i \iy\ \ \ {\rm as} \ \ t\to \pm\iy.
\]
The quasi-momentum $\vk$ maps the half plane
$\C_\pm=\{\l\in\C;\,\,\pm\Im\l>0\}$ onto the half-strip
$\K_\pm=\K\cap \C_\pm$ and $\s_{\rm
ac}(J^0)=\{\l\in\R;\,\,\Im\vk(\l)=0\}$.

Define the two strips $\K_S$ and $\cK$ by
$$
\K_S=-\K\qqq \mbox{and} \qqq \cK=\K_S\cup \K\ss
\{\vk\in \C: \Re \vk\in [-\pi,\pi]\}.
$$
The function $\vk$ has an analytic continuation from $\L_1\cap \C_+$
into $\L_1\cap \C_-$  through the infinite gaps $\g_q=(\l_q^-,\iy)$
by the symmetry and satisfies:

1)  $\vk$  is a  conformal mapping $\vk:\L_1\to \cK_+=\cK\cap \C_+$,
where we identify the boundaries $\{\vk=\pi+it, t>0\}$ and
$\{\vk=-\pi+it, t>0\}$.

2) $\vk:\L_2\to \cK_-=\cK\cap \C_-$  is a  conformal mapping, where
we identify the boundaries $\{\vk=\pi-it, t>0\}$ and $\{\vk=-\pi-it,
t>0\}$.

3) Thus $\vk:\L\to \cK$  is a  conformal mapping.

Consider the function $z=e^{i\vk(\l)},\,\, \l \in \L$. The function $z(\l),$ $\l \in \L,$
is a conformal mapping $z:\L\to \cZ=\C\sm \cup \ol g_j$, where the radial cut $g_j$ is given by
$$
g_j=(e^{-{h_j\/q}+i{\pi j\/q}}, e^{{h_j\/q}+i{\pi j\/q}}), \qqq j=\pm 1, ..., \pm (q-1).
$$

The function $z(\l),$  $\l \in \L,$ maps the first sheet $\L_1$ into the ``disk'' $\cZ_1=\cZ\cap \dD_1,$
$\dD_1=\{z\in\C:\,\,|z|<1\},$
and  $z(\cdot)$ maps the second sheet $\L_2$ into the domain  $\cZ_2=\cZ\sm \dD_1$.
In fact, we obtain the parametrization of the two-sheeted Riemann surface $\L$ by the ``plane''
$\cZ$. Thus below we call $\cZ_1$ also the ``physical sheet'' and $\cZ_2$ also the ``non-physical sheet''.

Note that if all $a_n^0=1, b_n^0=0$, then we have $\l={1\/2}(z+{1\/z})$. This function
$\l(z)$ is a conformal mapping from the disk $\dD_1$ onto the cut domain $\C\sm [-2,2]$.

Now, the functions $\psi^\pm(\l)$ can be considered as functions of $z\in \cZ$.
The functions $\psi^\pm_n(z)\equiv \psi^\pm_n(\l(z))$ are meromorphic in $\cZ$ with the only possible singularities at the images of the Dirichlet eigenvalues $z(\mu_j)\in\cZ$ and at $0.$ More precisely,\\
\no 1) $\psi_n^\pm$ are analytic in $\cZ \sm (\{z(\mu_j)\}_{j=1}^{q-1}\cup\{0\})$ and continuous up to $\pa \cZ\sm\{z(\mu_j)\}_{j=1}^{q-1}.$

\no 2) $\psi_n^\pm(z)$ has a simple pole at $z (\mu_j)\in\cZ$  if $\mu_j$ is a pole of $m_\pm,$ no pole if $\mu_j$ is not a singularity of $m_\pm$ (not a square root singularity if $\mu_j$ coincides with the band edge)  and if $\mu_j$  coincides with the band edge:  $\mu_j=\l_j^\sigma,$ $\sigma=+$ or $\sigma=-$, $j=1,\ldots, q-1,$ then
\[
\lb{B2psi} \psi_n^\pm(z)=\pm  \sigma(-1)^{q-j}\frac{iC(n)}{z -z(\l_j^\sigma)} +\cO(1),\qq\l\in [\l_{j-1}^+,\l_j^-],
\]
for some constant  $C(n)\in \R$. Note that the sign comes from the analytic continuation of the square root $\Omega(\l)$ using the definition
(\ref{branch}).

\no \no 3) The following identities hold true:
\[\lb{prop4}
\psi_n^\pm(\overline{z})=\psi_n^\pm(z^{-1})= \psi_n^\mp(z)=\overline{\psi_n^\pm(z)}\,\,\mbox{as}\,\, |z| =1.
\]
\no 4) The following asymptotics hold true:
$$
\psi_n^\pm(z)= (-1)^n\lt(\prod_{j=0}^{n-1}{}^*a_j\rt)^{\pm 1}z^{\pm
n}\lt(1 +\cO (z)\rt)
 \qqq \mbox{as} \qqq z\rightarrow 0.
 $$

We collect below some properties of the quasi-momentum $\vk$ on the gaps.

 On each $\g_j^+, j=0,1,\ldots,q,$ the
quasi-momentum $\vk(\l)$  has  constant real part and positive
$\Im\vk$:
 $$
 \Re \vk|_{\g_j^+}= -\frac{q-j}{q}\pi,
 \qqq\vk(\l_j^-)=\vk(\l_j^+)=-\frac{q-j}{q}\pi,\qqq
 \Im\vk|_{\g_j^+}>0.
 $$
  Moreover, as $\l$ increases from $\l_j^-$ to $\a_j$
  the  imaginary part  $\Im\vk\equiv h(\l)$
is monotonically increasing from $0$ to $h_j$ and as $\l$ increases
from $\a_j$ to $\l_j^-$   the  imaginary part  $\Im\vk\equiv h(\l+i0)$
is  monotonically  decreasing from $h_j$ to $0$. Then
 \[
 \lb{isin}
\frac12\vp_q(\l)(m_+(\l)-m_-(\l))= \sqrt{\D^2(\l)-1}=i\sin q\vk(\l)= -(-1)^{q-k}\sinh qh(\l+i0),
\]
where $\sinh qh=-2^{-1}(z^q-z^{-q})>0.$

\subsection{Scattering by finitely supported perturbations.}\lb{s-scattering}

For finitely-supported perturbation $(u,v)\in\cV_\nu$ we  define the Jost solutions  $f^\pm=(f^\pm_n)_{n\in \Z}$
 for the equation
\[
\lb{peq1}
a_{n-1} y_{n-1}+a_{n}y_{n+1}+b_ny_n=\l y_n,\qq  \,\, (\l,n)\in\C\ts\Z,
\]
(here $a_n=a_n^0+u_n,  b_n= b_n^0+v_n$ and $u_n=0, v_n=0$ for all $n\notin [0,p]$) by the conditions:
\[
\label{defJostline}
f_n^-=\p_n^-,\,\, \mbox{for}\,\, n \le 0\qqq \mbox{and} \qqq
f_n^+=\p_n^+,\,\,\mbox{for}\,\, n\ge p+1.
\]
Let $\vt^\pm, \vp^\pm $  be solutions to the equation \er{peq1}
satisfying the  conditions:
\[
\label{def2}
\vt_n^-= \vt_n, \ \ \vp^-_n=\vp_n\,\ \mbox{for}\,\,n\leq
0\qq\mbox{and}\qq  \vt_n^+=\vt_n, \ \vp_n^+=\vp_n\,\,\mbox{for}\,\,n\ge p+1.
\]
Note that each of $\vt_n^\pm, \vp_n^\pm, n \in \Z$ is a polynomial in $\l$.

The Jost solutions $f^\pm$  inherit the properties of $\p^\pm.$ We state this properties on the Riemann surface $\cZ$ as defined in Sections \ref{ss-cut}.
\begin{lemma}
\lb{l-jost} 1) Each $f_n^\pm, n\in\Z$, is analytic in $\cZ\sm \{0\}$
and continuous up to $\partial
\cZ\setminus\{z(\mu_j)\}_{j=1}^{q-1}.$ Moreover, the following
identities hold true:
\[
f^\s=\vt^\s +m_\s \vp^\s, \qqq \qqq \s=\pm.
\]
\[
 f_n^\pm(\overline{z})=f_n^\pm(z^{-1})=\overline{f_n^{\, \pm}(z)}\qq \mbox{for}\qqq |z| =1.
\]
2) $f_n^\pm(z)$ does not have a singularity at $z (\mu_j)$  if $\mu_j$ is not a singularity
(square root singularity if $\mu_j$ coincides with the band edge) of  $m_{\pm},$ otherwise, $f_n^\pm(z)$  can have either a simple pole at $z (\mu_j)$ if $\mu_j$ is a pole of $m_\pm,$ or a square root singularity,
\[
\lb{B2} f_n^\pm(\l)=\pm \sigma(-1)^{q-j}\frac{iC(n)}{\sqrt{\l -\l_j^\sigma}} +\cO(1),\qq \l\in [\l_{j-1}^+,\l_j^-],
\]
if $\mu_j$  coincides with the band edge:  $\mu_j=\l_j^\sigma,$ $\sigma=+$ or $\sigma=-$, $j=1,\ldots, q-1.$ Here
  $C(n)$ is bounded and real, the factor $ \sigma(-1)^{q-j}$ comes from the analytic continuation of the square root $\Omega(\l)$ using Definition
(\ref{branch}).
\end{lemma}

Define the unperturbed  Wronskian $\{\cdot,\cdot\}^0$ and the perturbed  Wronskian $\{\cdot,\cdot\}$
for sequences $f=(f_n)_{n\in \Z}, g=(g_n)_{n\in \Z}$ by
$$
\{f,g\}_n^0=a_n^0(f_ng_{n+1}-f_{n+1}g_n),\qqq
\{f,g\}_n=a_n(f_ng_{n+1}-f_{n+1}g_n).
$$
 Note that if  $f, g$ are solutions of \er{peq1}, then
the Wronskian $\{f,g\}_n$ is independent of $n$.

The Jost solutions  $f^\pm_n(\l)$ and
$\overline {f^\pm_n(\l)},$ $\l\in{\rm int}\,\sigma_{\rm ac}(H^0),$ are solutions of the same equation $Hf=\l
f$ and using $\overline{\psi_n^\pm(\l)}=\psi_n^\mp(\l),$ $\l\in\sigma_{\rm ac}(H^0),$ we have
\[
\lb{SM1}
\{f^\pm, \overline{f^\pm}\}=\{\psi^\pm,\psi^\mp\}^0=a_0^0(m_\mp -m_\pm)=
\mp a_0^0\frac{(z^q-z^{-q})}{\varphi_q},\,\,2i\Omega(\l)=z^q-z^{-q}.
\]

  We denote
  \begin{equation}\lb{def_of_sw}
  s=\{f^+, \overline{ f^-}\},\qqq  w=\{f^-, f^+\}.
\end{equation}

 Moreover we have (see \cite{IK1})
\begin{align}
& w=\{f_n^-,f_n^+\}=\const=\{f_0^-,f_0^+\}=a_0(f_0^-f_1^+-f_1^-f_0^+)=a_0f_1^++(v_0-a_0^0m_-)
 f_0^+,\lb{w1}\\
&  s=\{f_n^+,\overline{f_n^-}\}=\const=\{f_0^+,\overline{f_0^-}\}=
 a_0(f_0^+\overline{f_1^-}-f_1^+\overline{f_0^-})= (a_0^0\overline{m}_--v_0) f_0^+ -a_0f_1^+,\lb{s1}
\end{align}
 where we have used that $f_0^-=\psi_0^-=1$ and $a_0f_1^-=a_0^0m_--v_0,$ $\overline{m^\pm}=m^\mp,$ for $\l\in\sigma_{\rm ac}(H^0),$ by applying the Jacobi equation (\ref{pert}).

The following identities hold true:
\[
\lb{f-relationss}
f^\pm_n=\a \overline{f^\mp_n}+\b_\mp f^\mp_n,\qqq\l\in{\rm int}\,\sigma_{\rm ac}(H^0),
\]
where
\[
\lb{SM2}
\a ={\{f^\mp,f^\pm\}\/\{f^\mp,\overline{ f^\mp}\}}={\vp_q \{f^-, f^+\}\/a_0^0(z^q-z^{-q})}
={\vp_q w\/a_0^0(z^q-z^{-q})},
\]
\[
\lb{SM3}
\b_-={\{f^+,\overline{ f^-}\}\/\{f^-,\overline{ f^-}\}}={\vp_q s\/a_0^0(z^q-z^{-q})},
\qqq
\b_+={\{f^-,\overline {f^+}\}\/\{f^+,\overline{ f^+}\}}={\vp_q \overline {s}\/a_0^0(z^q-z^{-q})}
\]
since  $\overline{ s}=-\{f^-,\overline{ f^+}\}$.
Using (\ref{prop4}), we get $\overline{\b}_{\pm}=-\b_\mp$ for $\l\in{\rm int}\,\sigma_{\rm ac}(H^0)$
and
\begin{equation}
\lb{ws0}
|\alpha(\l)|^2 =1
+|\beta_\pm(\l)|^2,\,\,\l\in{\rm int}\,\sigma_{\rm ac}(H^0).
\end{equation}
Applying (\ref{SM2}), (\ref{SM3}) in  (\ref{ws0}) we get
\begin{equation}\lb{functional identity}
w(z)w(z^{-1})+\left(\frac{a_0}{\vp_q(z)}\right)^2(z^q-z^{-q})^2=s(z)s(z^{-1}),\qqq |z|=1,\qq z^2\neq 1.
\end{equation}

We define the scattering matrix
\begin{equation}\lb{Smatrix}
S(\l)=\left(
          \begin{array}{cc}
            T(\l) & R_-(\l) \\
            R_+(\l) & T(\l) \\
          \end{array}
        \right),\,\,\l\in\sigma (H^0),
\end{equation} for the pair $(H,H^0),$ where
$$T(\l)=\frac{1}{\alpha(\l)},\,\,R_\pm(\l)=\frac{\beta_\pm(\l)}{\alpha
(\l)}=\frac{\mp\{f^\mp(\l),\overline{ f^\pm(\l)}\}}{\{f^-(\l),
f^+(\l)\}}.$$
We have also $$R_+=\frac{\overline{s}}{w},\qq R_-=\frac{s}{w}.$$

 The matrix $S(\l)$ is unitary:
$|T(\l)|^2+|R_\pm(\l)|^2=1,$
$T(\l)\overline{R_+(\l)}=-\overline{T(\l)} R_-(\l),$ and extends to $\Lambda$ as a meromorphic function. The quantities
$T$ and $R_\pm$ are the transmission and the reflection coefficients
respectively:
$$T(\l)f^\pm_n(\l)=\left\{\begin{array}{lr}
                            T(\l) \psi^\pm_n(\l),& n\rightarrow \pm\infty \\
                            \psi_n^\pm(\l) +R_\mp(\l)\psi_n^\mp(\l),
                            & n\rightarrow \mp\infty
                          \end{array}\right.,\,\,\,\,\l\in\sigma
                          (H^0).
$$

The determinant of the scattering matrix is given by
$$\det S(\l)=T^2-R_+R_-=\frac{1}{\alpha^2} +\frac{|\beta|^2}{\alpha^2}=\frac{|\alpha|^2}{\alpha^2}=\frac{\overline{\alpha (\l)}}{\alpha (\l)}.$$

\section{Inverse resonance problem.}\lb{s-first}
\setcounter{equation}{0}

In this section  we prove the Theorems \ref{Th1} and \ref{zerosR+1}.

We consider the scattering data $(\hat{w},f_0^+).$ In Section \ref{s-pert} we summarize the characteristic properties of the scattering data and show (see
Lemma \ref{l_wf_in_class})  that, if  $V\equiv (u,v)\in\cV_\nu,$ then $(\hat{w},f_0^+)\in \gC_\nu.$

 We will show how from this data we can reconstruct the reflection coefficients and the norming constants which define a unique Jacobi operator $H.$ .

In Section \ref{s-invscat} we give an account of  the inverse scattering results from  \cite{Kh1} and \cite{EMT} in the form suitable for us.

In Section \ref{s-Th1} we prove Theorem \ref{Th1}. In Section \ref{s-zerosR+1} we give a sketch of the  proof of Theorem \ref{zerosR+1}.

\subsection{Characteristic properties of the scattering data.}\lb{s-pert}
\setcounter{equation}{0}
In this section we summarize the properties of the scattering data $(\hat{w},f_0^+)$ which are    needed for the proof of our main results.

Let $M_\pm\in\C$ denote (the projection of) the set of
poles  of $m_\pm.$ Let $M_{\rm e}$ denote the set of
 square root singularities of $m_\pm$ if $\mu_k=\l_j^\pm,$ $j=1,\ldots, q-1.$
 Note that $M_+\cap M_-=\emptyset.$
 We put
$$
\vp_q=a_0^0D^+D^-,\qqq D^+=\prod_{\mu_k\in M_+\cup M_{\rm e}} (\wt\l-\mu_k),\qqq
 D^-=\prod_{\mu_k\in M_-} (\wt\l-\mu_k),
$$ where $\,\wt {}\,:\,\,\Lambda\mapsto\C$ is the
natural projection introduced in (\ref{projection}). Note that this definition of $D^\pm$ differs from that used in \cite{IK1}.
We mark with $\hat{}$ the  modified (regularized) quantities:
$\hat{\psi}^\pm=D^\pm\psi^\pm,$ $\hat{f}^\pm =  D^\pm f^\pm,$ $\hat{w}=\frac{\vp_q}{a_0^0}w,$  which are analytic in $\Lambda_1.$ We denote also
$\hat{s}=\frac{\vp_q}{a_0^0}s,$ the meromorphic function on $\L_1.$ Note that the function $(D^+)^2s$ is analytic on $\L_1.$

In \cite{IK1} we proved  the following result.
\begin{lemma}\lb{l-3.2} Let $2i\Omega(\l)=2i\sin q\vk(\l)=z^q-z^{-q}$ and  $\l\in\sigma_{\rm ac}(H^0).$  The following identities hold true.
\begin{align}
& w=a_0f_1^++(v_0-a_0^0m_-)
 f_0^+,\,\,
  s=(a_0^0m_+-v_0) f_0^+ -a_0f_1^+,\qq \hat{w}=\frac{\vp_q}{a_0^0}w,\qq\hat{s}=\frac{\vp_q}{a_0^0}s,\nonumber\\
  &\hat{w}\hat{w}^* =4\Omega^2+\hat{s}\hat{s}^*,\qq w+s=\frac{2i\Omega}{\vp_q}f_0^+,\,\,  \hat{w}+\hat{s}=\frac{2i\Omega}{a_0^0}f_0^+,\lb{ws}\\
&\hat{w}=2i\Omega (1+A)-J,\qq \hat{s}=-2i\Omega (1+\tilde{A})+\tilde{J},\nonumber\\
&
 A= \frac12\left(\frac{a_0}{a_0^0}\vp_1^++\frac{v_0}{a_0^0}\vp_0^++\vt_0^+\right)-1=\frac12\left[
\left(\frac{ a_0}{a_0^0}\vp^+_1-\vp_1\right)+\frac{ v_0}{a_0^0}\vp_0^++(\vt^+_0-\vt_0)
\right],\nonumber\\
& \tilde{A}= \frac12\left(\frac{a_0}{a_0^0}\vp_1^++\frac{v_0}{a_0^0}\vp_0^+-\vt_0^+-\frac{2\phi}{\vp_q}\vp_0^+\right)-1,\nonumber\\
&J=-\left[\frac{ a_0}{a_0^0}\vp_q\vt^+_1+\phi(\frac{ a_0}{a_0^0}\vp^+_1-\vt^+_0)+
\frac{ v_0}{a_0^0}(\vp_q\vt^+_0+\phi
\vp^+_0)+\vt_{q+1}\vp^+_0 \right],\,\,\vt_{q+1}\vp_0^+=-\frac{\phi^2+\Omega^2}{\vp_q}\vp^+_0,\nonumber\\
&\tilde{J}=-\left[\frac{ a_0}{a_0^0}\vp_q\vt^+_1+\phi(\frac{ a_0}{a_0^0}\vp^+_1-\vt^+_0)+
\frac{ v_0}{a_0^0}(\vp_q\vt^+_0+\phi
\vp^+_0)-\frac{\phi^2-\Omega^2}{\vp_q}\vp^+_0 \right],\nonumber\\
&A-\tilde{A}=\vt_0^++\frac{\phi}{\vp_q}\vp_0^+,\qqq J-\tilde{J}=\frac{2\Omega^2}{\vp_q}\vp_0^+,\nonumber\\
&4\Omega^2(1+A)^2+J^2=4\Omega^2+4\Omega^2(1+\tilde{A})^2+\tilde{J}^2 \lb{analyticcont}.
\end{align}
\end{lemma}
{\bf Remark.} From Lemma \ref{l-3.2} it follows that $A,J$ are polynomials, whereas functions $\tilde{A},\tilde{J}$ are rational with simple poles at $\mu_n\neq \l_j^\pm,$ $n=1,\ldots,q-1,$ $j=1,\ldots, q,$ where $\mu_n$ is a zero of $\vp_q$ and $\l_j^\pm$ is  an endpoint of a finite gap.

Note also that $\vp_j,\vt_j,$ $j=1,2$ are related via  the wronskian property: $a_0(\vt_0^+\vp_1^+-\vt_1^+\vp_0^+)=a_0^0(\vt_0\vp_1-\vt_1\vp_0)=a_0^0.$

In order to define the class of scattering data we need the  properties summarized in  the following lemma.

\begin{lemma}
\lb{L1.2} Suppose $(u,v)\in \cV_\n$, where $\n\in \{2p,2p-1\}$ and  $\{\rho_k\}_{k=1}^N=\s_{\rm bs}(H)$ is the set of bound states of $H.$ Let $\{\mu_j\}_{j=1}^{q-1}$ be the set of zeros of $\vp_q.$ Then  for all $k=1,\ldots,N$ the functions $\hat{f}^\pm,$ $\hat{w}$ and $\hat{s}$ satisfy the following properties:\\
\begin{equation}
\lb{1.16}
\phantom{=}\hspace{-1cm}\mbox{ 1)}\,\,\hat{f}^\pm_n(\rho_k)=\beta_\mp(\rho_k) \frac{D^\pm(\rho_k)}{D^\mp(\rho_k)}\hat{f}^\mp_n(\rho_k),\qq \hat{s}(\rho_k)\overline{\hat{s}(\rho_k)}=-4(1-\Delta^2(\rho_k)).
\end{equation}
2) Denote $\displaystyle{ c^\pm_k=\frac{\hat{f}^\pm(\rho_k)}{\hat{f}^\mp(\rho_k)},\,\,\gamma_{\pm,k}=\left(\sum_{n\in\Z}|\hat{f}_n^\pm(\rho_k)|^2\right)^{-1}}.$ Then
\begin{equation}\lb{6.11}
c^\pm_k=\beta_\mp \frac{D^\pm}{D^\mp}(\rho_k),\qq\gamma_{\pm,k}=
-\left(c^\pm_k\hat{w}'(\rho_k)\right)^{-1}>0,\qq \gamma_{+,k}\gamma_{-,k}=\left(\hat{w}'(\rho_k)\right)^{-2}.\end{equation}\\
3) For $k=1,\ldots,N,$ if   $\rho_k\in\gamma_j^+$ for some $j=0,\ldots,q,$ then  \begin{equation}\lb{gamma_plus}
 \gamma_{+,k}=-\frac{D^-}{D^+}\frac{2i\Omega(\rho_k)}{\hat{s}\hat{w}'(\rho_k)}=\frac{D^-}{D^+}\frac{2(-1)^{q-j}\sinh qh(\rho_k)}{\hat{s}\hat{w}'(\rho_k)}>0,
 \end{equation} where $h(\rho_k)=\Im\vk(\rho_k)>0.$\\
\\
4) If $\mu_j\neq \l_j^\pm$ for all $j=1,\ldots q-1,$  then we have
\begin{equation}
\lb{nakoncah} \hat{s}(\l_j^\pm)=\overline{\hat{s}(\l_j^\pm)}=\tilde{J}(\l_j^\pm)=-\hat{w}(\l_j^\pm)=J(\l_j^\pm);
\end{equation}
if $\mu_j=\l_j^\sigma,$ for some $j=1,\ldots,q-1$ and $\sigma=+$ or $\sigma=-,$  then \begin{equation}\lb{nakoncah1}\hat{s}(\l_j^\sigma)=\tilde{J}(\l_j^\sigma)=\hat{w}(\l_j^\sigma)=-J(\l_j^\sigma).
\end{equation}\\
5) The function
$\hat{w}$ is real on $\R$ and
\begin{equation}
\lb{bound_below}
|\hat{w}(z)|\geq |z^q -z^{-q}|\qq\mbox{for any}\qq |z|=1.
\end{equation}
\end{lemma}
{\bf Proof.} 1) Relations in (\ref{1.16}) follow from the analytic continuation of the identities (\ref{f-relationss})  and (\ref{ws0})
(see also (\ref{analyticcont}))
  as $\alpha (\rho_j)=0.$

2) Formulas in (\ref{6.11}) come from (6.11) in \cite{EMT}.

3)  Formula (\ref{gamma_plus}) follows using that
$\beta_-=\frac{D^+D^-}{(z^q-z^{-q})}s(z)=\frac{1}{(z^q-z^{-q})}\frac{D^-}{D^+}\hat{s}(z)$ and (\ref{isin}).

4) Proof of (\ref{nakoncah}).\\
Suppose $\mu_j\neq\l_j^\pm$ for all $ j=1,\ldots q-1.$ Then $m_-(\l_{j}^\pm)= \overline{m_-(\l_{j}^\pm)}$ and (\ref{nakoncah})  follows from  (\ref{w1}), (\ref{s1}).  Moreover,  as $f^\pm(\l_j^\pm)$ are real, it follows also that $\overline{\hat{s}(\l_{j}^\pm)}=-\hat{w}(\l_{j}^\pm).$

 Now, suppose $\mu_j=\l_j^-,$ for some $j=1,\ldots,q-1.$   Then
   $f_n^\pm(\l)$ are pure imaginary in the limit $\l\rightarrow \l_j^--$ (see (\ref{B2})) and we
  use the formulas $w(\l)=\{f_n^-,f_n^+\},$ $s(\l)=\{f_n^+,\overline{{f}_n^-}\}.$ Taking the  limit $\l\rightarrow \l_j^--$ we get
 $\hat{s}(\l_j^-)=\hat{w}(\l_j^-).$ The case $\mu_j=\l_j^+$ follows similarly. \\
5) The property (\ref{bound_below}) $|\hat{w}(z)|\geq |z^q -z^{-q}|$ for $|z|=1$  follows from $\overline{z}=z^{-1}$ and $\overline{(z^q-z^{-q})}=-(z^q-z^{-q})$ as in (\ref{ws}) we have
$(z^q-z^{-q})^2=-|z^q-z^{-q}|^2.$

\hfill\qed

Using Lemma \ref{L1.2} and the asymptotics  of the Jost functions given  in \cite{IK1}   we have the following result.
\begin{lemma}\lb{l_wf_in_class} If $(u,v)\in\cV_\nu,$ $\nu\in\{2p-1,2p\},$ then $(\hat{w},f_0^+)\in\gC_\nu$ and
\begin{align*}
&   \hat{w}=2i\Omega(1+A)-J=\ca \frac{A_p}{c_1}\l^q\left(1+{\mathcal O}(\l^{-1})\right) \qq &{\rm if}
\qq   \l\in\L_1\\
 -\frac{v_0}{A_p} c_2\l^{\nu+q-1}\left(1+{\mathcal O}(\l^{-1})\right) \qq &{\rm
if}\qq \l\in\L_2 \ac \qqq{\rm as}\qq \l\to \iy,\\
&f_0^+=\vt_0^++\frac{\phi}{\vp_q}\vp_0^++ i\frac{\Omega(\l)}{\vp_q}\vp_0^+=\ca c_1A_p+{\cO}(\l^{-1}) \qq &{\rm if}
\qq   \l\in\L_1\\
 -\frac{c_2 }{A_p}\l^{\nu}+\cO(\l^{\nu-1}) \qq &{\rm
if}\qq \l\in\L_2 \ac \qqq{\rm as}\qq \l\to \iy,
\end{align*}
with the constants $c_1,c_2$ given in (\ref{c2_2p-1}) and $A_p$  defined in (\ref{defAp}).

The polynomials $1+A,$ $J$ have asymptotics
\begin{align*}
& 1+A=- \frac{c_2}{2A_p}v_0\l^{\nu-1}\left(1+{\mathcal O}(\l^{-1})\right),\qqq J= \frac{c_2}{2A_p}v_0\l^{\nu+q-1}\left(1+{\mathcal O}(\l^{-1})\right).
\end{align*}

\end{lemma}

At last we reformulate some properties of the functions $\hat{w}$ and $\hat{s}$ on the Riemann surface $\cZ.$
\begin{lemma}
\lb{l-asymptotics}
Let $(u,v)\in\cV_\nu,$ $\nu\in\{2p-1,2p\}.$ Then the function $z^q\hat{w}(z),$ is entire   in $\cZ,$  the function $z^{q-1}\hat{s}(z)$ is analytic in  $\cZ_1$ and have poles in $\cZ_2.$ The functions $\hat{w},$ $\hat{s}$ satisfy  the functional equation
\begin{equation}\lb{ws}
 \hat{w}(z)\hat{w}(z^{-1})+(z^q-z^{-q})^2=\hat{s}(z)\hat{s}(z^{-1}),\qqq |z|=1,\,\,z^2\neq 1.
\end{equation}
Moreover, the following asymptotics hold
\begin{align*}
&\hat{w}=\frac{A_p}{c_1}z^{-q}\left[1 +{\mathcal O}\left( z\right)\right],\qqq \hat{s}=-\frac{A_p}{c_1}v_0z^{1-q}\left[ 1 +{\mathcal O}\left(z\right)\right]\qq \mbox{as}\qq z\to 0,\\
&\hat{w}=-\frac{c_2}{A_p}v_0 z^{\nu+q-1}\left[ 1 +{\mathcal O}\left(z^{-1}\right)\right],\qqq \hat{s}=\frac{c_2}{A_p} z^{\nu+q}\left[ 1 +{\mathcal O}\left(z^{-1}\right)\right]\qq \mbox{as}\qq z\to \infty,
\end{align*}
with the constants $c_1,c_2$ given in (\ref{c2_2p-1}) and $A_p$ is defined in (\ref{defAp}).
\end{lemma}
{\bf Proof.} Identity (\ref{ws}) follows from (\ref{functional identity}) using that for $|z|=1$ we have $\vp_q(z^{-1})=\vp(\overline{z})=\vp(z).$
The asymptotics of $\hat{w}$ follows from Lemma \ref{l_wf_in_class} using $2\Delta=z^q+z^{-q}=\l^q+{\mathcal O}(\l^{q-1})$ as $\lambda\rightarrow\infty.$
  Similarly follows the asymptotics for $\hat{s}.$
\qed

\subsection{Inverse scattering problem}\lb{s-invscat}
 We consider the relation between the {\bf left/right scattering data} $S_\pm(H)$ for $H,$
$$S_\pm(H)=\{ R_\pm (z),\,\, z\in \S^1;\,\,\rho_k, \gamma_{\pm,k}>0,\,\,
k=1,\ldots,N\},$$ and the perturbation coefficients $(u,v)$ in the Jacobi operator $H.$

In \cite{Kh1} and \cite{EMT} the inverse scattering problem was solved for Jacobi operators which are short range perturbations of  periodic (quasi-periodic in \cite{EMT}) finite-gap operators. Here
we give a short summary of their results   in the context of periodic background {\bf with finitely supported perturbations.}  In this case  the proofs follows straightforward from the methods  in \cite{EMT} and we omit them.

{\bf First we consider the direct problem.}
Let $\S^1$ denote the unit circle $|z|=1$ and consider the measure on $\S^1$
\begin{equation}\lb{measure}
d\omega(z)=\prod_{j=1}^{q-1}\frac{\l(z)-\mu_j}{\l(z)-\alpha_j}\frac{dz}{z},
 \end{equation}
 where $\alpha_j\in\gamma_j$ is the zero of $\Delta'(\l)$ (see Section \ref{ss-cut} and \cite{EMT}).

 Introduce the transformation operator $K_\pm$ by
$$
(K_\pm h)_n=\sum_{m=n}^{\pm\infty}K_\pm(n,m)h_m,\qq h=(h_n)_{n\in \Z},
$$
where the kernel $K_\pm(n,m)$ is given by
\begin{equation}\lb{5.1}
K_\pm(n,m)=\frac{1}{2\pi i}\int_{|z|=1}f^\pm_n(z)\psi^\mp_md\omega(z).
\end{equation}
The kernels $K_\pm(n,m),$ $n,m\in\Z,$ are the Fourier coefficients of the Jost solution $f^\pm_n$ with respect to the orthonormal system $\{\psi_n^\pm\}_{n\in\Z}$ in the Hilbert space $L^2(\S^1,\frac{1}{2\pi i}d\omega).$

\begin{lemma}\lb{l5.1}
Assume $(u,v)\in\cV_\nu.$ Then the Jost solutions  $f^\pm$ have the form
\begin{equation}
\lb{5.4}
f^\pm_n=\sum_{m=n}^{\pm\infty}K_\pm(n,m)\psi_m^\pm,\qq |z|=1,
\end{equation}
where the kernels $K_\pm (n,m)$ of the finite rank operator  $K_\pm$ satisfy
$$
K_\pm(n,m)=0, \,\,\mbox{for}\,\,\pm m <\pm n,
$$
and
\begin{equation}
\lb{5.5}
|K_+(n,m)|\leq C\sum_{j=\max\{M,0\}}^pQ_j,\qq|K_-(n,m)|\le C\sum_{j=0}^{\min\{M,p\}}Q_j,\qq\mbox{for}\,\,\pm m >\pm n,
\end{equation}
where $Q_j=|u_j|+|v_j|$ and $M=[{n+m\/2}]+1,$ and the constant $C\equiv C(H^0)$ depends on the unperturbed periodic operator $H^0.$
\end{lemma}

\begin{lemma}
\label{l5.3}
Assume $(u,v)\in\cV_\nu.$ Then $a_n, b_n, n\in\Z,$ satisfy
\[
\frac{a_n}{a_n^0}=\frac{K_+(n+1,n+1)}{K_+(n,n)}=\frac{K_-(n,n)}{K_-(n+1,n+1)},
\lb{5.27}
\]
\[
 v_n=a_n^0\frac{K_+(n,n+1)}{K_+(n,n)}-a_{n-1}^0\frac{K_+(n-1,n)}{K_+(n-1,n-1)}
=a_{n-1}^0\frac{K_-(n,n-1)}{K_-(n,n)}-a_{n}^0\frac{K_-(n+1,n)}{K_-(n+1,n+1)}.\nonumber
\]
\end{lemma}

Let
\[
F^\pm(l,m)=F_0^\pm(l,m)+\sum_{j=1}^N\gamma_{\pm,j}\hat{\psi}^\pm_l(\rho_j)
\hat{\psi}_m^\pm(\rho_j),
\lb{7.7}
\]
\[
 {F_0}^\pm(l,m)=\frac{1}{2\pi i}\int_{|z|=1}
R_\pm(z) \psi^\pm_l(z) \psi_m^\pm(z)d\omega(z).
\lb{7.3}
\]
Note that $F_0^\pm(l,m)=F_0^\pm(m,l)$ is real.
The function $K_\pm(n,m)$ satisfies the equation
\begin{equation}\lb{T(7.5)}
K_\pm(n,m)+\sum_{l=n}^{\pm\infty}K_\pm(n,l){F_0}^\pm(l,m)=\frac{\delta(n,m)}{K_\pm(n,n)}-
\sum_{j=1}^N\gamma_{\pm,j}\hat{f}^\pm_n(\rho_j)\hat{\psi}^\pm_m(\rho_j).
\end{equation}
We define the Gel'fand-Levitan-Marchenko operator
$$
(\cF_n^\pm f)(j)=\sum_{l=0}^\iy F^\pm(n\pm l,n\pm j)f_l,\qqq f=(f_l)_0^\iy\in\ell^\infty(0,\infty).
$$
\begin{theorem}\lb{Th7.1}
The kernel $K_\pm(n,m)$ of the transformation operator satisfies  the Gel'fand-Levitan-Marchenko equation
\begin{equation}\lb{7.9}
(1+\cF_n^\pm)K_\pm(n,n\pm .)=(K_\pm(n,n))^{-1}\delta_0,
\end{equation}
where $K_\pm(n,n)=\langle\delta_0,(1+ \cF_n^\pm)^{-1}\delta_0\rangle^{\frac12},$
and we have
 \begin{equation}\lb{7.10}
|F^+(n,m)|\leq C\sum_{j=\left[\frac{n+m}{2}\right]+1}^{p}\left(|u_j|+|v_j|\right),
\qq|F^-(n,m)|\leq C\sum_{j=0}^{\left[\frac{n+m}{2}\right]-
1}\left(|u_j|+|v_j|\right),
\end{equation}
where the constant $C=C(H^0)$ is of the same nature as in (\ref{5.5}).
\end{theorem}

Now we recall the procedure which allows the reconstruction of the perturbation coefficients $(u,v)$ for the Jacobi operator $H$ from the {\bf left/right scattering data} $S_\pm(H)$ for $H,$
$$
S_\pm(H)=\{ R_\pm (z),\,\,|z|=1;\,\,\rho_k, \gamma_{\pm,k}>0,\,\,
k=1,\ldots,N\}.
$$ We summarize the properties of the scattering data $S_\pm(H)$ for $H$ with the finitely supported perturbation coefficients $(u,v)\in\cV_\nu.$\\
{\bf Hypothesis 1}
{\it
 The scattering data $S_\pm(H)$ satisfy
the following conditions:\\
 (i) The reflection coefficients $R_\pm(z)$ are continuous
except possibly at $z_l=z(E_l),$ where
$\{E_l\}_{l=0}^{2q-1}\equiv\{\l^\pm_k\}_{k=0}^q$ and fulfill
$\displaystyle \overline{R_\pm(z)}=R_\pm(\overline{z}).$
Moreover, $\displaystyle |R_\pm(z)|<1$ for $z\neq z_l$
and
$$1-|R_\pm(z)|^2\geq C\prod_{l=0}^{2q-1}|z -z_l|^2.$$\\
The functions $F_\pm(n,m)$ satisfy
\begin{equation}\lb{7.10}
F_\pm(n,m)=0\,\, \mbox{for}\,\,\pm(n+m)\geq M.
\end{equation}
 for some $M\in\Z.$\\
(ii) The values $\rho_k\in\R\setminus\sigma_{\rm ac}\,(H^0),$ $1\leq
k\leq N,$ are distinct and the norming constants $\gamma_{\pm,k},$
$1\leq k\leq N,$ are positive.\\
(iii) $T(z)$ defined via Poisson-Jensen type formula ((6.25) in
\cite{EMT}) extends to a single values function on $\cZ_1$
(i.e. it has equal
values on the corresponding slits).\\
(iv) Transmission and reflection coefficients satisfy
$$
\lim_{z\,\,\rightarrow\,\,z_l}\sqrt{\Delta^2(z)-1}\,\frac{R_\pm(z)+1}{T(z)}=0,\qq
z_l\neq z (\mu_j),\,\,j=1,\ldots,q-1,
$$$$
\lim_{z\,\,\rightarrow\,\,z_l}\sqrt{\Delta^2(z)-1}\,\frac{R_\pm(z)-1}{T(z)}=0,\qq
z_l=z (\mu_j),\,\,j=1,\ldots,q-1,
$$ where $z_l=z(E_l),$
and the consistency conditions
$$
\frac{R_-(z)}{R_+(\overline{z})}=-\frac{T(z)}{T(\overline{z})},\qqq
\gamma_{+,k}\gamma_{-,k}=\frac{({\rm
Res}\,_{\rho_k}T(\l))^2}{4(\Delta^2(\rho_k)-1)}.
$$
}

\begin{theorem}\lb{Th7.5}
Suppose that Hypothesis 1 is satisfied. Then, for $n\in\Z,$ the Gel'fand-Levitan-Marchenko operator $\cF_n^\pm
:\,\,\ell^2\to\ell^2$ has finite rang. Moreover, $1+\cF_n^\pm$
is positive and hence invertible.

In particular, the Gel'fand-Levitan-Marchenko equation (\ref{7.9}) has a unique solution  and
 $S_+(H)$ or $S_-(H)$ uniquely determines $H$ and the  finitely supported perturbation
 $(u,v).$
\end{theorem}

{\bf Inverse problem.} If $S_\pm$ (satisfying Hypothesis 1 {\em
(i),(ii)}) and $H^0$ are known, we can construct $F^\pm (l,m)$ via
formula (\ref{7.7}) and thus derive the Gel'fand-Levitan-Gel'fand-Levitan-Marchenko
equation, which has a unique solution by Theorem \ref{Th7.5}. This solution
$$K_\pm(n,n)=\langle\delta_0,(1+ \cF_n^\pm)^{-1}\delta_0\rangle^{\frac12},\qq K_\pm(n,n\pm j)=\frac{1}{K_\pm(n,n)}\langle\delta_j,(1+ \cF_n^\pm)^{-1}\delta_0\rangle^{\frac12}$$
is the kernel of the transformation operator. Since $1+\cF_n^\pm$ is positive, $K_\pm(n,n)$ is positive and we can set in accordance with Lemma \ref{l5.3}
\[
a_n^+=a_n^0\frac{K_+(n+1,n+1)}{K_+(n,n)},\qq b_n^+=b_n +
a_n^0\frac{K_+(n,n+1)}{K_+(n,n)}-a_{n-1}^0\frac{K_+(n-1,n)}{K_+(n-1,n-1)},
\lb{uv1}
\]
\[
a_n^-=a_n^0\frac{K_-(n,n)}{K_-(n+1,n+1)},\qq b_n^-=b_n +
a_{n-1}^0\frac{K_-(n,n-1)}{K_-(n,n)}-a_{n}^0\frac{K_-(n+1,n)}{K_-(n+1,n+1)}.
\lb{uv2}
\]
Let $H^+,$ $H^-$ be the associated Jacobi operators.

\begin{lemma}\label{l8.1}
Suppose that a given set $S_\pm$ satisfies Hypothesis 1 (i)-(ii). Then the sequences $(a_n^\pm-a_n^0,b_n^\pm-b_n^0)_{n\in\Z},$  defined in (\ref{uv1}), (\ref{uv2}) have finite support.\\  Moreover,  $f^\pm_n=\sum_{m=n}^{\pm\infty}K_\pm(n,m)\psi_m^\pm,$ where $K_\pm(n,m)$ is the solution of the Gel'fand-Levitan-Marchenko equation (\ref{7.9}), satisfies $Hf^\pm=\l f^\pm$ and $f_n^\pm=\psi_n^\pm$ for $\pm n > n^\pm,$ for some $n^\pm\in\N.$
\end{lemma}

Now in \cite{EMT} it is shown that $a_n^+=a_n^-$ and $b_n^+=b_n^-$ and we have
\begin{theorem}\lb{Th8.2} Hypothesis 1 is necessary and
sufficient for a set $S_\pm$ to be left/right scattering data of a
unique Jacobi operator $H$ associated with sequences $a,b$ such that $(u,v)\in\cV_\nu.$
\end{theorem}

We set \begin{equation}\lb{uv}
\frac{a_n^0+u_n}{a_n^0}=\frac{K_+(n+1,n+1)}{K_+(n,n)},\qq v_n=
a_n^0\frac{K_+(n,n+1)}{K_+(n,n)}-a_{n-1}^0\frac{K_+(n-1,n)}{K_+(n-1,n-1)},
\end{equation}

\subsection{Proof of Theorem \ref{Th1}.}\lb{s-Th1}
   Let $(u,v)\in\cV_\nu.$ Then the scattering data satisfy Hypothesis 1.  By Lemma \ref{L1.2} we have also  $(\hat{w},f_0^+)\in\gC_\nu.$    Lemma \ref{L1.2} yields the norming constants $\gamma_{\pm,k},$ $k=1,\ldots,N.$
The following lemma shows the inverse relation.
\begin{lemma}\lb{l_hyp} Suppose $(\hat{w},f_0^+)\in\gC_\nu,$ $\nu\in\{2p-1,2p\},$ and $(u,v)$ be defined by (\ref{uv}). Then\\
1) Hypothesis 1 is satisfied.\\
2) $(u,v)\in\cV_\nu.$
Moreover,  $f^\pm_n=\sum_{m=n}^{\pm\infty}K_\pm(n,m)\psi_m^\pm,$ where $K_\pm(n,m)$ is the solution of the Gel'fand-Levitan-Marchenko equation (\ref{7.9}), satisfies $Hf^\pm=\l f^\pm$ and $f_n^-=\psi_n^-$ for $n\leq 0,$
$f_n^+=\psi_n^+$ for $n\geq p+1.$

\end{lemma}
{\bf Proof.}\\
If $(\hat{w},f_0^+)\in\gC_\nu,$   then  we have $$R_+=\frac{s(z^{-1})}{w(z)}=\frac{\hat{s}(z^{-1})}{\hat{w}(z)},\qq |z|=1.$$
1) We need to check that the scattering data $S_\pm(H)$ satisfies  Hypothesis 1.\\
(i) As $R_-=s/w,$ then we have
$$
1-|R_-|^2=\frac{a_0^0|z^q-z^{-q}|^2}{\vp_q^2|w|^2}\equiv |T|^2=\frac{|z^q-z^{-q}|^2}{|\hat{w}|^2}\geq C\prod_{l=1}^{2q-1}|z-z(E_l)|^2.
$$
It follows from the fact that $\sup\hat{w} <\const$ for $|z|=1.$

We prove  (\ref{7.10}). If $(u,v)\in\cV_\nu$ then  for $l+m\geq 2p$ we have
$$
{F_0}^+(l,m)=-\sum_{j=1}^N\gamma_{+,j}\hat{\psi}_l^+(\rho_j)
\hat{\psi}_m^+(\rho_j)
$$
(see proof of 2) below) and
\begin{equation*}
F^+(l,m)={F_0}^+(l,m)+\sum_{j=1}^N\gamma_{+,j}\hat{\psi}^+_l(\rho_j)
\hat{\psi}_m^+(\rho_j)=0,\qq l+m\geq 2p,
\end{equation*}
 and similar for $F^-$.\\
(ii) follows from 3) in the definition of $\gC_\nu.$ \\
(iii) follows as the transmission coefficient $T$ is defined via $T(z)=\frac{z^q-z^{-q}}{\hat{w}(z)}.$\\
 (iv) We have \cite{EMT}
$$
\sqrt{\Delta^2(z)-1}\,\frac{R_+(z)+1}{T(z)}=
{w(z)+s(z^{-1})\/2}\prod_{j=1}^{q-1}(\l-\mu_j)
$$$$
\sqrt{\Delta^2(z)-1}\,\frac{R_-(z)+1}{T(z)}=
{(w+s)\/2}\prod_{j=1}^{q-1}(\l-\mu_j).
$$

If $\mu_j\neq E_l$ then $f^\pm$ are continuous and real at $\l=E_l$ and the two Wronskians cancel (see \cite{EMT}). This also follows from Lemma \ref{L1.2}, equation (\ref{nakoncah}).

Otherwise, if $\mu_j = E_l,$ the Wronskians are purely imaginary (by property (\ref{B2})) and add up. This also follows from   Lemma \ref{L1.2},
equation (\ref{nakoncah1}).

We have for $|z|=1$
$$
\frac{R_-(z)}{R_+(\overline{z})}=\frac{s(z)w(\overline{z})}{w(z))
\overline{s(\overline{z})}},\qq \frac{T(z)}{T(\overline{z})}=-\frac{\vp_q(\overline{z})w(\overline{z})}{\vp_q(z)w(z)},
$$
which together with $\overline{s(z)}=s(\overline{z}),$  $\overline{w(z)}=s(\overline{w})$  and $\vp_q(\overline{z})=\vp_q(z)$ for $|z|=1$ give the first consistency condition.

Now we have
$$
4(\Delta^2-1)=(z^q-z^{-q})^2,\,\, T(z)=\frac{z^q-z^{-q}}{\hat{w}(z)}\,\Rightarrow\,\,T^2=\frac{4(\Delta^2-1)}{\hat{w}^2},
$$$$
\hat{w}(\l)=
\hat{w}'(\rho_j)(\l-\rho_j)+{\mathcal O}(\l-\rho_j)^2\qq\Rightarrow\qq ({\rm Res}_{\rho_j}T)^2=\frac{4(\Delta^2(\rho_j)-1)}{(\hat{w}'(\rho_j))^2}.
$$
As from Lemma \ref{L1.2}, equation (\ref{6.11}), we have $$\gamma_{+,j}\gamma_{-,j}=\left(\hat{w}'(\rho_j)\right)^{-2},$$ then we have the second consistency condition.

2) Recall  Formula (\ref{7.3})
 $$
F_0^+(l,m)=\frac{1}{2\pi i}\int_{|z|=1}
R_+(z) \psi^+_l(z) \psi_m^+(z)d\omega(z),
$$
where $d\omega(z)$ is given in (\ref{measure}).

Observe that $d\omega$ is meromorphic on $\cZ_1$ with simple pole at $z=0.$ In particular, there are no poles at $z(\alpha_j).$
To evaluate the integral we use the residue theorem. Take a closed contour in $\cZ_1$ and let this contour approach $\partial \cZ_1.$ The function $R_\pm(z) \psi^\pm_l(z) \psi_m^\pm(z)$ is continuous on $\{|z|=1\}\setminus\{z(E_j)\}$ and meromorphic on $\cZ_1$ with simple poles at $z(\rho_j)$ and eventually a pole at $z=0.$

Due to the properties of $\hat{s},$ $\hat{w},$ $\psi_l^+$ we have
$$ R_+\sim_{z\sim 0}=\frac{z^{-(\nu+q)}}{z^{-q}}=z^{-\nu},\qq\psi_l^+\psi_m^+\sim_{z\sim 0}z^{l+m}.$$ Suppose $l+m\geq \nu+1$ ($+1$ is due to singularity $z^{-1}$ in $d\omega$). Then the integrand is bounded near $z=0$ and we
apply the residue theorem to the only poles at the eigenvalues.

We have (\cite{EMT}, (3.23)) $$\frac{dz}{d\l}=z\frac{\prod_{j=1}^{q-1}(\l-\alpha_j)}{2(\Delta^2(\l)-1)^{1/2}}$$ and if $z_j=z(\rho_j)$ then
${\rm Res}_{z=z_j} F(z)=z'(\rho_j){\rm Res}_{\l=\rho_j}F(z(\l)).$

Then  we get
$$
F_0^\pm(l,m)=\sum_{j=1}^N{\rm Res}_{\rho_j}\left(\frac{R_\pm(z)D^+D^-\psi_l^\pm(\l)\psi_m^\pm(\l)}{2(\Delta^2(\l)-1)^{1/2}}\right)=
\sum_{j=1}^N{\rm Res}_{\rho_j}\left(\frac{D^\mp}{D^\pm}\frac{R_\pm(z)\hat{\psi}_l^\pm(\l)\hat{\psi}_m^\pm(\l)}
{2(\Delta^2(\l)-1)^{1/2}}\right),
$$
where we used $\prod_{j=1}^{q-1}(\l-\mu_j)=D^+D^-.$
Now we consider ${F_0}^+$ only, as calculations for $F_0^-$ are similar. We have
$$
F_0^+(l,m)=\sum_{j=1}^N{\rm Res}_{\rho_j}\left(\frac{D^-}{D^+}\frac{\hat{s}(z^{-1})\hat{\psi}_l^+(\l)
\hat{\psi}_m^+(\l)}{\hat{w}(z)(z^q-z^{-q})}\right).
$$
The functions $z^{\nu+q}\hat{s}(z^{-1}),$ $z^q\hat{w}(z)$ are analytic.  Thus
 $$
 z^{\nu+q}\hat{s}(z^{-1})=z_j^{\nu+q}\hat{s}(z_j^{-1})+{\mathcal O}(z-z_j),\qq
 z^q\hat{w}(z)=z_j^q\hat{w}'(\rho_j)(\l-\rho_j)+{\mathcal O}(\l-\rho_j)^2.
 $$
The function $\frac{\hat{\psi}_l^+(\l)
\hat{\psi}_m^+(\l)}{(z^q-z^{-q})}$ is bounded at $\rho=\rho_j$ and
  we write
 $$
 \frac{\hat{s}(z^{-1})\hat{\psi}_l^+(\l)
\hat{\psi}_m^+(\l)}{\hat{w}(z)(z^q-z^{-q})}=\frac{z^{\nu+q}\hat{s}(z^{-1})
}{z^q\hat{w}(z)}\cdot \frac{\hat{\psi}_l^+(\l)
\hat{\psi}_m^+(\l)}{z^q-z^{-q}}\cdot z^{-(\nu+q)+q}
$$
  We have $\hat{\psi}_l^\pm(z)\sim z^{\pm l}$ as $z\rightarrow 0.$  Now if $l+m\geq \nu+1$ (see above, $+1$ comes from $d\omega$) then $\hat{\psi}_l^+(\l)
\hat{\psi}_m^+(\l)z^{-\nu-1}$ is bounded. Using that
 $\hat{s}(z_j^{-1})=(z_j^q-z_j^{-q})^2(\hat{s}(z_j))^{-1}$ (due to Lemma \ref{L1.2}),
 we get
$$
F_0^+(l,m)=\sum_{j=1}^N\frac{D^-(\rho_j)}{D^+(\rho_j)}
\frac{(z_j^q-z_j^{-q})\hat{\psi}_l^+(\rho_j)
\hat{\psi}_m^+(\rho_j)}{\hat{w}'(\rho_j)\hat{s}(z_j)}=-\sum_{j=1}^N
\gamma_{+,j}\hat{\psi}_l^+(\rho_j)
\hat{\psi}_m^+(\rho_j).
$$
Then equation (\ref{7.7}) implies
\begin{equation*}
F^+(l,m)=F_0^+(l,m)+\sum_{j=1}^N\gamma_{+,j}\hat{\psi}^+_l(\rho_j)
\hat{\psi}_m^+(\rho_j)=0,\qq l+m\geq \nu +1,
\end{equation*}
and the Gel'fand-Levitan-Marchenko equation
$$
K_+(n,m)+\sum_{l=n}^{+\infty} K_+ (n,l) F^+
(l,m)=\frac{\delta_{nm}}{K_+(n,n)},\qq  m\geq n,
$$
implies that,  if $n+m\geq \nu +1,$ the  kernel of the transformation operator
$K_\pm(n,m)$ satisfies
$$
K_+(n,m) =\frac{\delta_{nm}}{K_+(n,n)},\qq  m\geq n,\qq m+n\geq \nu +1.
$$
Thus we get the following  properties:
$$\mbox{if}\,\,2n\geq \nu +1,\,\,\mbox{then}\,\,K_0(n,n)=\pm 1;\qq \mbox{if}\,\,n+m\geq \nu+1,\,\,m\neq n,\,\,\mbox{then}\,\,K_0(n,m)=0.$$

As by (5.27) in \cite{EMT} we have  $$\frac{a_n}{a_n^0}=\frac{K_+(n+1,n+1)}{K_+(n,n)},\qq v_n=
a_n^0\frac{K_+(n,n+1)}{K_+(n,n)}-a_{n-1}^0\frac{K_+(n-1,n)}{K_+(n-1,n-1)},$$ and if $\nu=2p,$
then, $a_n=a_n^0$ for $n\geq p+1,$  and $v_n=0$ for $n\geq p+1.$ The case $\nu=2p-1$ is similar.

Analogously it follows  that $u_n,$ $v_n$ are $0$ for $n<0.$ \qed

{\bf Proof of Theorem \ref{Th1}.  Uniqueness.}
If $V\equiv(u,v)\in\cV_\nu,$ then Lemma \ref{l_wf_in_class}  implies that $(\hat{w}, f_0^+)\in\gC_\nu,$ which yields a mapping $V\mapsto (f_0^+,\hat{w})$ from $\cV_\nu$ into $\gC_\nu.$

We will show  uniqueness. Let $V\equiv(u,v)\in\cV_\nu.$ Then Lemma \ref{l_wf_in_class} gives unique $(\hat{w}, f_0^+)\in\gC_\nu.$    Lemma \ref{L1.2} yields the norming constants $\gamma_{\pm,j},$ $j=1,\ldots,N.$ By 1) in Lemma \ref{l_hyp} the Hypothesis 1 is satisfied. Then, by Theorem \ref{Th8.2}, these data determine the finitely supported perturbation uniquely. Then we deduce that the mapping $V\mapsto (\hat{w}, f_0^+)$ is an injection.

{\bf Surjection} of the mapping $V\equiv(u,v)\mapsto (\hat{w}, f_0^+)$ follows from 2) in Lemma \ref{l_hyp}.
\subsection{ Proof of Theorem \ref{zerosR+1}.}\lb{s-zerosR+1} First we note that, as the zeros of $\hat{w}$  and $R_-+1$ are disjoint from the end-points $\l_k^\pm,$ $k=0,\ldots,q,$  and from the points $\mu_j,$ $j=1,\ldots,q-1,$ such that $\vp_q (\wt\mu_j)=0$, then the zeros of $$R_-+1=\frac{w+s}{w}=\frac{2i\Omega f_0^+}{\vp_q w}=\frac{2i\Omega f_0^+}{a_0^0\hat{w}}$$ (see (\ref{ws})) coincide with the zeros of $f_0^+.$   We used also that the zeros of $f_0^+$ and $w$ are disjoint, which follows from the Wronskian property: as
$w=\{f^-,f^+\}=a_0(f_0^-f_1^+-f_1^-f_0^+),$ then $w(\l_0)=0,$ $f_0^+(\l_0)=0$ would imply $f_0^-(\l_0)=0,$ which is impossible.

Moreover, the multiplicities of zeros of $R_-+1$ and $f_0^+$ coincide, i.e. the zeros are simple.

Now using the zeros of $f_0^+,$ the polynomial $\vp_0^+$ and the constant $c_3$ we can reconstruct the polynomial $F(\l)=\vp_qf_0^+(f_0^+)^*,$ $\l\in\L_1.$ We recall that  $f^*(\l)=\overline{f(\overline{\l})}.$
The zeros of $f_0^+$ and the polynomial $F$ were studied in detail in \cite{IK3}. Applying Theorem 1.4 in \cite{IK3} we know that   the function $f_0^+$ is uniquely determined by the polynomials $F$ and $\vp_0^+,$ supposing that the zeros of $f_0^+$ are simple.

In \cite{IK1} we show that the set of zeros of the polynomial $\cF=\hat{w}\hat{w}^*$ on $\C$ coincide with the projection of the set of states of $H:$ $\wt\sigma_{\rm st}(H),$ with the same multiplicities. Thus, as the polynomial $F$ in \cite{IK3}, the polynomial $\cF$ can be reconstructed from $\wt\sigma_{\rm st}(H)$ and the constant $c_3v_0$ (see (\ref{tn})). Now applying the polynomial interpolation method (see \cite{A})
as in \cite{IK3} we can reconstruct $\hat{w}=2i\Omega(1+A)-J$ from $\sigma_{\rm st}(H)$ using the polynomials $A$ and $\cF.$ Now we have also the function
$\hat{s}=\frac{2i\Omega}{a_0^0}f_0^+-\hat{w},$ see (\ref{ws}).\qed

\end{document}